\documentstyle[12pt,amscd,amssymb]{amsart}

\newcommand{\sk}{\vspace{1mm}}
\renewcommand{\backslash}{/}

\newcommand{\ezlk}{E_{\z}^{l_{\z}-k,k}}
\newcommand{\lz}{l_{\z}}
\newcommand{\Mm}{{\frak M}_H(c_1,c_2)}
\newcommand{\MM}{{\frak M}}

\newcommand{\wc}{\d_{S,\z}^{w,d}}

\newcommand{\exseq}[3]{0 \ar #1 \ar #2 \ar #3 \ar 0}
\newcommand{\seq}[3]{{#1}_{#2}, \ldots, {#1}_{#3}}

\newcommand{\inc}{\hookrightarrow}
\newcommand{\ar}{\rightarrow}

\newcommand{\x}{\times}
\newcommand{\iso}{\cong}
\newcommand{\isom}{\stackrel{\sim}{\ar}}

\renewcommand{\pt}{\text{pt}}
\newcommand{\vol}{\text{vol}}
\newcommand{\CP}{{\Bbb C \Bbb P}}

\newcommand{\End}{\text{End}}

\newcommand{\Hilb}{\text{Hilb}}
\newcommand{\Jac}{\text{Jac}}

\newcommand{\rk}{\text{rk}}

\newcommand{\Hom}{\text{Hom}}
\newcommand{\Sym}{\text{Sym}}

\newcommand{\PD}{\text{P.D.}}

\newcommand{\ch}{\text{ch}\:}
\newcommand{\Ext}{\text{Ext}}

\newcommand{\cC}{{\cal C}}

\newcommand{\cE}{{\cal E}}
\newcommand{\cF}{{\cal F}}

\newcommand{\cH}{{\cal H}}
\newcommand{\cM}{{\cal M}}
\newcommand{\cN}{{\cal N}}
\newcommand{\cL}{{\cal L}}
\newcommand{\cO}{{\cal O}}

\newcommand{\cU}{{\cal U}}

\newcommand{\cZ}{{\cal Z}}
\renewcommand{\AA}{{\Bbb A}}

\newcommand{\HH}{{\Bbb H}}
\newcommand{\PP}{{\Bbb P}}
\newcommand{\QQ}{{\Bbb Q}}
\newcommand{\RR}{{\Bbb R}}
\newcommand{\SS}{{\Bbb S}}
\newcommand{\ZZ}{{\Bbb Z}}

\renewcommand{\a}{\alpha}
\renewcommand{\b}{\beta}
\renewcommand{\d}{\delta}
\newcommand{\g}{\gamma}
\newcommand{\e}{\varepsilon}
\newcommand{\f}{\epsilon}

\renewcommand{\l}{\lambda}

\newcommand{\s}{\sigma}

\renewcommand{\o}{\omega}

\newcommand{\r}{\rho}

\newcommand{\z}{\zeta}

\renewcommand{\S}{\Sigma}
\newcommand{\D}{\Delta}
\renewcommand{\L}{\Lambda}

\newcommand{\fra}{{\frak a}}

\newcommand{\frM}{{\frak M}}

\theoremstyle{plain}
\begingroup
\theorembodyfont{\sl}
\newtheorem{thm}{Theorem}
\newtheorem{cor}[thm]{Corollary}
\newtheorem{lem}[thm]{Lemma}
\newtheorem{prop}[thm]{Proposition}
\endgroup

\theoremstyle{definition}

\theoremstyle{remark}
\newtheorem{rem}[thm]{Remark}

\newtheorem{notation}[thm]{Notation}

\title[Wall-crossing formulae]{Wall-crossing formulae 
for algebraic surfaces with $q >0$}

\author{Vicente Mu\~noz}
\address{Departamento de \'Albegra, Geometr\'{\i}a y Topolog\'{\i}a  \\
Facultad de Ciencias \\
Universidad de M\'alaga \\ Campus de Teatinos, s/m \\ 29071 M\'alaga \\ Spain}
\email{vmunoz@@agt.cie.uma.es}
\date{August, 1997 }

\begin{document}

\maketitle

\begin{abstract}
  We extend the ideas of Friedman and Qin~\cite{flips} to find the
  wall-crossing
  formulae for the Donaldson invariants of
  algebraic surfaces with $p_g=0$, $q>0$ and anticanonical divisor 
  $-K$ effective, for any
  wall $\z$ with $\lz={1 \over 4} (\z^2-p_1)$ being $0$ or $1$. 
\end{abstract}

\section{Introduction}
\label{sec:intro}

The Donaldson invariants of a smooth oriented $4$-manifold $X$ 
depend by definition on a Riemannian metric $g$. In the
case $b^+>1$ they however turn out to be independent of $g$. When
$b^+=1$, they depend on $g$ through a structure of walls and chambers, that we 
recall briefly here (we refer to~\cite{Kotschick1}~\cite{KM2} for more details).

Fix $w \in H^2(X;\ZZ)$. Then for any $p_1 \leq 0$ with $p_1 \equiv w^2
\pmod 4$, we set $d=-p_1-{3 \over 2}(1-b_1 +b^+)$, for half of the 
dimension of the moduli space $\cM_{X,g}^{w,d}$ of $g$-antiselfdual connections
on the $SO(3)$-principle bundle over $X$ 
with second Stiefel-Whitney class the reduction 
mod $2$ of $w$, and first Pontrjagin number $p_1$. The corresponding 
Donaldson invariant will be denoted $D_{X,g}^{w,d}$. 
This is a linear functional on the 
elements of degree $2d$ of $\AA(X)=\Sym^*(H_0(X) \oplus H_2(X)) \otimes
\bigwedge^*(H_1(X)\oplus H_3(X))$, where the degree of elements
in $H_i(X)$ is $4-i$ ($H_i(X)$ will always denote homology with rational 
coefficients, and similarly for $H^i(X)$). 
This invariant is only defined in principle for generic metrics.

From now on let $X$ be a compact smooth oriented $4$-manifold with $b^+=1$. 
Let $\HH$ be the image of the positive cone $\{ x \in H^2(X; \RR)/ x^2>0 \}$ in 
$\PP(H^2(X; \RR))$, which is a model of the hyperbolic disc of dimension
$b^-$. The period point of $g$ is the line $\o_g \in \HH \subset 
\PP(H^2(X; \RR))$ given by the selfdual harmonic forms for $g$.
A {\bf wall} of type $(w,p_1)$ is a non-empty hyperplane
$W_{\z}=\{ x \in \HH / x\cdot \z =0 \}$ in $\HH$,
with $\z \in H^2(X; \ZZ)$, such that $\z \equiv w\pmod 2$ 
and $p_1 \leq \z^2 < 0$. The connected components of
the complement of the walls of type $(w,p_1)$ in $\HH$ are the 
{\bf chambers} of type $(w,p_1)$. 

Let $\frM$ denote the space of metrics of $X$. Then we have a map $\frM \to 
\HH$ which sends every metric $g$ to its period point $\o_g$. The connected 
components of the preimage of the chambers of $\HH$ are, by definition, the chambers
of $\frM$. A wall $W'_{\z}$ for $\frM$ is a non-empty preimage of a wall $W_{\z}$ 
for $\HH$.
When $g$ moves in a chamber $\cC'$ of $\frM$ 
the Donaldson invariants do not change. But when
it crosses a wall they change. 
So for any chamber $\cC'$ of $\frM$, we have defined
$D_{X,\cC'}^{w,d}$, by choosing any generic metric $g \in \cC'$, 
so that the moduli space is 
smooth, and computing the corresponding Donaldson invariants 
(to avoid flat connections we might have to
use the trick in~\cite{MM}). For a path of metrics $\{ g_t
\}_{t \in [-1,1]}$, with $g_{\pm 1} \in \cC'_{\pm}$, we have the 
difference term $\d_{X}^{w,d}(\cC'_-,\cC'_+)=D_{X,\cC'_+}^{w,d}-D_{X,\cC'_-}^{w,d}$.

When $b_1=0$, Kotschick and 
Morgan~\cite{KM2} prove that the invariants only depend on the chamber $\cC$ 
of $\HH$ in which the period point of the metric lies. For this, they find that the 
change in the Donaldson invariant when the metric crosses a wall $W'_{\z}$
depends only on the class
$\z$ and not on the particular metric 
having the reducible antiselfdual connection (Leness~\cite{Leness} 
points out that their argument 
is not complete and checks that it is true at least for the case 
$\lz={1 \over 4} (\z^2-p_1) \leq 2$). 
In this case, the difference term is defined as
$$
  \d_{X}^{w,d}(\cC_-,\cC_+)=D_{X,\cC_+}^{w,d}-D_{X,\cC_-}^{w,d},
$$
for chambers $\cC_{\pm}$ of $\HH$. Then $\d_{X}^{w,d}(\cC_-,\cC_+)=
\sum \wc$, where the sum is taken over all $\z$ defining walls separating
$\cC_-$ and $\cC_+$.
Moreover $\wc=\e(\z,w) \d_{S,\z}^d$, with $\d^d_{S,\z}$ not dependent on $w$,
$\e(\z,w)=(-1)^{({\z-w \over 2})^2}$.

Now suppose $S$ is a smooth algebraic surface (not necessarily with 
$b_1=0$), endowed with a Hodge metric $h$
corresponding to a polarisation $H$. 
Let $\Mm$ be the Gieseker compactification of the 
moduli space of $H$-stable rank two bundles $V$ on
$X$ with $c_1(V)=\cO(L)$ (a fixed line bundle with topological first
Chern class equal to $w$) and $c_2={1 \over 4}(c_1^2 -p_1)$. 
The Donaldson invariants (for the metric $h$) can be computed using $\Mm$
(see~\cite{FM}) whenever the moduli spaces $\Mm$ are generic
(i.e. $H^0(\End_0E)=H^2(\End_0E)=0$, for every stable bundle $E \in \Mm$).
The period point of $h$ is the line spanned by $H \in H^2(X;\ZZ) \subset H^2(X;\RR)$.
Now let $C_S \subset \HH$  be the image of the {\bf ample cone} of $S$, 
i.e. the subcone of the positive cone
generated by the ample classes (polarisations). We have walls and
chambers in $C_S$ in the same vein as before (actually they are the intersections
of the walls and chambers of $\HH$ with $C_S$, 
whenever this intersection is non-empty). 
Now $\Mm$ is constant on the chambers of $C_S$
(and so the invariant stays the same), and when $H$ crosses
a wall $W_{\z}$, $\Mm$ changes (see~\cite{Qin}).
From the point of view of the Donaldson invariants, this corresponds to 
restricting our attention from the positive cone of $S$ to its ample cone.

When the irregularity $q$ of $S$ is zero, the wall-crossing terms have been
found out in~\cite{flips}~\cite{Gottsche-notengo}~\cite{Gottsche-bott}. 
In~\cite{flips}
Friedman and Qin obtain some wall-crossing formulae for algebraic
surfaces $S$ with $-K$ being effective ($K=K_S$ the canonical
divisor) and the irregularity $q=0$
(equivalently, $b_1=0$). 
We want to adapt their results to the case
$q >0$ modifying their arguments where necessary. 
If $-K$ is effective then the change of $\Mm$ when $H$ crosses 
a wall $W$ can 
be described by a number of flips. We shall write
the change of the Donaldson invariant as a sum of
contributions $\wc$, for the different $\z$ defining $W$.

\begin{rem}
\label{rem:kld}
The condition of $-K$ being effective can be relaxed for the case $q=0$ to
the following two conditions: $K$ is not effective, $\pm 
\z +K$ are not effective for any $\z$ defining the given wall (we call such
a wall a {\bf good wall}, see~\cite{Gottsche-bott}~\cite{Gottsche-notengo}).
Probably the same is true for the case $q>0$, since these two conditions ensure
that the change in $\Mm$ when crossing a wall is described by flips. Nonetheless
we will suppose $-K$ effective, which allows us to define the Donaldson invariants
for any polarisation. Note that when $-K$ is effective, all walls are good.
\end{rem}

The paper is organised as follows. In section~\ref{sec:2}
we extend the arguments of~\cite{flips} to the case $q>0$. 
In sections~\ref{sec:lz=0} and~\ref{sec:lz=1} we compute 
the wall crossing formulae for any wall with
$l_{\z}={1 \over 4}(\z^2-p_1)$ being $0$ and $1$ respectively.
Then in section~\ref{sec:5},
we give the two leading terms of the wall crossing 
difference for any wall $\z$. As a consequence of our results,
we propose a conjecture on the shape of 
the wall crossing terms.
In the appendix we give, for the convenience of the reader, a list of all 
the algebraic surfaces with $p_g=0$ and 
$-K$ effective, i.e. the surfaces to which the results from this
paper apply.

\noindent {\em Acknowledgements:\/} I am very grateful to my D. Phil.\
supervisor Simon Donaldson, for many good ideas. Conversations with Lothar 
G\"ottsche have been very useful for a checking of the formulae here obtained.
Also I am indebted to the Mathematics 
Department in Universidad de M\'alaga for their hospitatility and financial 
support.

\noindent {\em Note:\/} After the completion of this work, L.\ G\"ottsche
provided me with a copy of~\cite{Gottsche-notengo}. The arguments for
computing the wall-crossing terms
in~\cite{Gottsche-notengo} can also be extended to the case $q>0$, in a
similar fashion to the work carried out in this paper.

\section{Wall-crossing formulae}
\label{sec:2}

From now on, $S$ is a smooth algebraic manifold with irregularity $q \geq 0$ 
and $p_g=0$ (equivalently $b^+=1$) and with anticanonical divisor 
$-K$ effective. Let $w \in H^2(S;\ZZ)$, $p_1 \equiv w^2 \pmod 4$. 
Put 
$$
  d=-p_1- {3\over 2}(1-b_1+b^+)=-p_1-3(1 -q)
$$ 
and let $\z$ define a wall of type $(w,p_1)$. In every
chamber $\cC$ of the ample cone,
we have well-defined the Donaldson invariant $D^{w,d}_{S,\cC}$
associated to polarisations in that chamber.
For two different chambers $\cC_+$ and
$\cC_-$, there is a {\bf wall-crossing difference term\/}
$$
   \d_S^{w,d}(\cC_-,\cC_+)=D^{w,d}_{S,\cC_+}- D^{w,d}_{S,\cC_-},
$$
which can be written as a sum
$$
   \d_S^{w,d}(\cC_-,\cC_+)=\sum_{\z} \d_{S,\z}^{w,d},
$$
where $\z$ runs over all walls of type $(w, p_1)$ with $\cC_- \cdot \z <0 
< \cC_+ \cdot \z$.

Suppose from now on that $\cC_-$ and $\cC_+$ are two adjacent chambers
separated by a single wall $W_{\z}$ of type $(w,p_1)$. 
For simplicity, we will assume that the wall $W_{\z}$ is only
represented by the pair
$\pm \z$ since in the general case we only need to add
up the contributions for every pair representing the wall.
Then the wall-crossing term is $\wc$.
Set 
$$
  l_{\z}=(\z^2 -p_1)/4 \in \ZZ.
$$

Let $\z$ define the wall separating $\cC_-$ from $\cC_+$ and put, 
as in~\cite[section 2]{flips}, $E_{\z}^{n_1,n_2}$ to be the set of all
isomorphism classes of non-split extensions of the form
$$
  \exseq{\cO(F) \otimes I_{Z_1}}{V}{\cO(L-F) \otimes I_{Z_2}},
$$
where $F$ is a divisor such that $2F-L$ is homologically
equivalent to $\z$, 
and $Z_1$ and $Z_2$ are two zero-dimensional subschemes of $S$ with
$l(Z_i)=n_i$ and such that $n_1+n_2=l_{\z}$.
Let us construct $E_{\z}^{n_1,n_2}$ explicitly. Consider $H_i
=\text{Hilb}_{n_i}(S)$ the Hilbert scheme of $n_i$ points on $S$, 
$J=\text{Jac}^F(S)$ the Jacobian parametrising
divisors homologically equivalent to $F$, $\cZ_i \subset S \x H_i$
the universal codimension $2$ scheme, and $\cF \subset S \x J$ the universal
divisor. Then we define $\cE_{\z}^{n_1,n_2} \to J\x H_1 \x H_2$ to be
$$
  \cE=\cE_{\z}^{n_1,n_2}=\cE{\text{xt}}_{\pi_2}^1(\cO_{S \x (J\x H_1
  \x H_2)} (\pi_1^* L -\cF) \otimes I_{\cZ_2}, \cO_{S \x (J\x H_1 \x
  H_2)} (\cF) \otimes I_{\cZ_1}),
$$ 
for $\pi_1:S \x (J\x H_1 \x H_2) \ar S$ and $\pi_2:S \x 
(J\x H_1 \x H_2) \ar J\x H_1 \x H_2$, the projections (we do not
denote all pull-backs of sheaves explicitly).
This is a vector bundle over $J\x H_1 \x H_2$ of rank 
$$
  \rk(\cE)=l_{\z}+h^1(\cO_S(2F-L))=l_{\z}+h(\z)+q,
$$ 
where 
$$
  h(\z)= {\z \cdot K_S \over 2} -{\z^2 \over 2} -1,
$$
by Riemann-Roch~\cite[lemma 2.6]{flips}. Note that $l_{\z} \geq
0$ and $h(\z)+q \geq 0$. Put
$N_{\z}=\rk(\cE)-1$. Then $E_{\z}^{n_1,n_2}
=\PP((\cE_{\z}^{n_1,n_2})^{\vee})$ (we
  follow the convention $\PP(\cE)= \text{Proj} (\oplus_i
  S^i(\cE))$),
which is of dimension $q+2l_{\z}+(l_{\z}+h(\z)+q)$. 
Also $N_{\z}+N_{-\z}+q+2l_{\z}=d-1$. We will have to
treat the case $\rk(\cE)=0$ (i.e. $l_{\z}=0$ and $h(\z)+q=0$) separately.

We can modify the arguments in sections 3 and 4 of~\cite{flips}
to get intermediate moduli
spaces $\MM_0^{(k)}$ together with embeddings $E_{\z}^{l_{\z}-k,k} \inc
\MM_0^{(k)}$ and $E_{-\z}^{k,l_{\z}-k} \inc
\MM_0^{(k-1)}$, fitting in the following diagram
$$
\begin{array}{ccccccccccccc}
  && \widetilde{\MM}_0^{(l_{\z})} & & & & \cdots & & & &
    \widetilde{\MM}_0^{(0)} & &\\ 
  & \swarrow & & \searrow && \swarrow & &\searrow & & \swarrow & &
    \searrow & \\  
  \MM_0^{(l_{\z})} & & & & \MM_0^{(l_{\z}-1)} & & & & \MM_0^{(0)} & &
    & & \MM_0^{(-1)} \\ 
  \parallel &&&&&&&&&&&& \parallel \\
  \MM_- &&&&&&&&&&&& \MM_+
\end{array}
$$
where $\widetilde{\MM}_0^{(k)} \ar \MM_0^{(k)}$ is the blow-up of
$\MM_0^{(k)}$ at $E_{\z}^{l_{\z}-k,k}$ and $\widetilde{\MM}_0^{(k)}
\ar \MM_0^{(k-1)}$ is the blow-up of $\MM_0^{(k-1)}$ at
$E_{-\z}^{k,l_{\z}-k}$. This is what is called a flip.
Basically, the space $E_{\z}= \sqcup E_{\z}^{l_{\z}-k,k}$ parametrises
$H_-$-stable sheaves 
which are $H_+$-unstable. Analogously, $E_{-\z}= \sqcup
E_{-\z}^{k,l_{\z}-k}$ parametrises $H_+$-stable sheaves 
which are $H_-$-unstable. Hence one could say that $\MM_+$ is obtained from
$\MM_-$ by removing $E_{\z}$ and then attaching $E_{-\z}$. The picture above
is a nice description of this fact and allows us the find the universal
sheaf for $\MM_+$ out of the universal sheaf for $\MM_-$ by a sequence
of elementary transforms.

The point is that whenever $-K_S$ is effective,
we have an embedding $E_{\z}^{0,\lz} \ar \MM_-$ 
(the part of $E_{\z}$ consisting of bundles)
and rational maps $E_{\z}^{k,\lz-k} \dashrightarrow 
\MM_-$, $k>0$, but if we blow-up
$\MM_-$ at $E_{\z}^{0,\lz}$, we have already an embedding from 
$E_{\z}^{1,\lz-1}$ to this latter space. Now we can proceed inductively
for $k=0,\ldots, \lz$. Analogously, we can have started from $\MM_+$ 
blowing-up $E_{-\z}^{k,l_{\z}-k}$ one by one. The diagram above says
that we can perform 
these blow-ups and blow-downs alternatively, instead of
first blowing-up $\lz+1$ times and then blowing-down $\lz+1$ times.
We see that the exceptional divisor in $\widetilde{\MM}_0^{(k)}$ is a 
$\PP^{N_{\z}} \x \PP^{N_{-\z}}$-bundle 
over $J \x H_{\lz-k} \x H_{k}$.

When adapting the arguments of~\cite[sections 3 and 4]{flips},
the only place requiring serious changes
is proposition 3.7 in order to prove proposition 3.6.

\begin{prop}{\em {\bf (\cite[proposition 3.6]{flips})}}
\label{prop:lod}
  The map $E^{\lz-k,k}_{\z} \ar \MM_0^{(\z, {\bf k})}$ is an immersion. The
  normal bundle $\cN_{\z}^{l_{\z}-k,k}$ to $E^{\lz-k,k}_{\z}$ in 
  $\MM_0^{(\z, {\bf k})}$ is exactly $\rho^*\cE^{k,l_{\z}-k}_{-\z} \otimes
  \cO_{E^{\lz-k,k}_{\z}}(-1)$, where $\rho: E^{\lz-k,k}_{\z} \ar J \x
  H_{\lz-k} \x H_k$ is the projection. Here we have defined
  $\cE_{-\z}^{k,\lz-k}=\cE{\text{xt}}_{\pi_2}^1(\cO_{S \x (J\x H_1 \x H_2)} 
 (\cF) \otimes I_{\cZ_1},
  \cO_{S \x (J\x H_1 \x H_2)} (\pi_1^* L -\cF) \otimes I_{\cZ_2})$. 
\end{prop}

Proposition~\ref{prop:lod} is proved as~\cite[proposition 3.6]{flips} making
use of the following
analogue of~\cite[proposition 3.7]{flips}

\begin{prop}
  For all nonzero $\xi \in \Ext^1=\Ext^1(\cO(L-F) \otimes 
  I_{Z_2},\cO(F) \otimes I_{Z_1})$, the natural map from a neighbourhood 
  of $\xi$ in $E^{\lz-k,k}_{\z}$ to $\MM_0^{(\z, {\bf k})}$
  is an immersion at $\xi$. The image of $T_{\xi} 
  E^{\lz-k,k}_{\z}$ in $\Ext^1_0(V,V)$ (the tangent space to 
  $\MM_0^{(\z, {\bf k})}$ at $\xi$, where $V$ is the sheaf corresponding to
  $\xi$)
  is exactly the kernel of the natural map $\Ext_0(V,V) \ar \Ext^1(\cO(F)  
  \otimes I_{Z_1},\cO(L-F) \otimes I_{Z_2})$, and the normal space to 
  $E^{\lz-k,k}_{\z}$ at $\xi$ in $\MM_0^{(\z, {\bf k})}$ 
  may be canonically identified with $\Ext^1(\cO(F)  
  \otimes I_{Z_1},\cO(L-F) \otimes I_{Z_2})$.
\end{prop}

\begin{pf}
We have
that $\Ext^1 (I_Z, I_Z)$ parametrises infinitesimal deformations of
$I_Z$ as a sheaf. The deformations of $I_Z$ 
are of the form $I_{Z'} \otimes \cO(D)$ for $D
\equiv 0$. The universal space parametrising these sheaves is
$\Hilb_r(S) \x \Jac^0(S)$, where $r$ is the length of $Z$. 
There
is an exact sequence 
$$
  \exseq{H^0 (\cE \text{xt}^1(I_Z,I_Z))}{\Ext^1 (I_Z, I_Z)}{H^1(\cH
   \text{om}(I_Z,I_Z))},
$$
where $H^0 (\cE \text{xt}^1(I_Z,I_Z))=H^0 (\cH\text{om}(I_Z,\cO_Z))=
\Hom (I_Z,\cO_Z)$ is the tangent space to $\Hilb_r(S)$ and 
$H^1(\cH \text{om}(I_Z,I_Z))=H^1(\cO)$ is the tangent space to the
Jacobian.
Analogously, $\Ext^1(V,V)$ is the space of infinitesimal deformations
of $V$ (but the determinant is not preserved). The infinitesimal
deformations preserving the determinant are given by the kernel
$\Ext_0^1(V,V)$ of a map $\Ext^1(V,V) \ar H^1(\cH\text{om}(V,V)) \ar H^1(\cO)$.
Now $E=E_{\z}^{l_{\z}-k,k}$ sits inside the bigger space $\tilde{E}=
\tilde{E}_{\z}^{l_{\z}-k,k}$ given as
$$
  \PP(\cE{\text{xt}}_{\pi_2}^1(\cO_{S \x (J_1\x H_1
  \x J_2 \x H_2)} (\pi_1^* L -\cF_2) \otimes I_{\cZ_2}, \cO_{S \x (J_1\x
  H_1 \x J_2\x H_2)} (\cF_1) \otimes I_{\cZ_1})^{\vee}),
$$
for $J_1=J_2=J$, $\cF_i \subset S \x J_i$ the universal divisor, and
$H_i$ the Hilbert scheme parametrising $Z_i$.
The arguments in~\cite[proposition 3.7]{flips} go through to prove
that for every non-zero $\xi \in \Ext^1=\Ext^1(\cO(L -F) \otimes
I_{Z_2}, \cO(F) \otimes I_{Z_1})$
we have the following commutative diagram with exact rows and columns
$$
\begin{CD}
  T_{\xi} E @>>> \Ext^1_0(V,V)  @>>> \Ext^1(\cO(F) \otimes
   I_{Z_1},\cO(L -F) \otimes I_{Z_2}) \\
   @VVV       @VVV      @| \\
  T_{\xi} \tilde{E}  @>>> \Ext^1(V,V)  @>>>\Ext^1(\cO(F) \otimes
       I_{Z_1},\cO(L -F) \otimes I_{Z_2}) \\ 
    @VVV         @VVV   @. \\
   H^1(\cO) @= H^1(\cO) @. 
\end{CD}
$$
So the natural map from a neighbourhood of $\xi$ in $E$ to
$\MM_0^{(\z, {\bf k})}$ is
an immersion at $\xi$ and the normal space may be canonically
identified with $\Ext^1(\cO(F) \otimes 
I_{Z_1},\cO(L -F) \otimes I_{Z_2})$.
\end{pf}

Therefore proposition~\ref{prop:lod} is true for $q>0$. The set up is now
in all ways analogous to that of~\cite{flips}.
We fix some notations~\cite[section 5]{flips}:

\begin{notation}
\label{not:wall}
Let $\z$ define a wall of type $(w,p_1)$.
\begin{itemize}
\setlength{\itemsep}{0pt}
  \item $\l_k$ is the tautological line bundle over $E_{\z}^{\lz-k,k}
  =\PP((\cE_{\z}^{\lz-k,k})^{\vee})$.
  $\l_k$ will also be 
  used to denote its first Chern class.
  \item $\r_k:S \x \ezlk \ar S \x (J \x H_{\lz-k} \x H_k)$ is the
  natural projection.
  \item $p_k: \widetilde{\MM}_0^{(k)} \ar  \MM_0^{(k)}$ is the blow-up
  of $\MM_0^{(k)}$ at $\ezlk$.
  \item $q_{k-1}: \widetilde{\MM}_0^{(k)} \ar  \MM_0^{(k-1)}$ is the
  contraction of $\widetilde{\MM}_0^{(k)}$ to $\MM_0^{(k-1)}$.
  \item The normal bundle of $\ezlk$ in $\MM_0^{(k)}$ is $\cN_k =
  \r_k^* \cE_{-\z}^{k,l_{\z}-k} \otimes \l_k^{-1}$, where
  $\cE_{-\z}^{k,\lz-k}=\cE{\text{xt}}_{\pi_2}^1(\cO_{S \x (J\x H_1 \x H_2)} 
 (\cF) \otimes I_{\cZ_1},
  \cO_{S \x (J\x H_1 \x H_2)} (\pi_1^* L -\cF) \otimes I_{\cZ_2})$. 
  \item $D_k=\PP(\cN_k^{\vee})$ is the exceptional divisor in
  $\widetilde{\MM}_0^{(k)}$.
  \item $\xi_k=\cO_{\widetilde{\MM}_0^{(k)}}(-D_k)|_{D_k}$ is the
  tautological line bundle on $D_k$.
  \item $\mu^{(k)}(\a)=-{1 \over 4}p_1({\frak g}_{\cU^{(k)}}) 
  \backslash \a$, for
  $\a \in H_2(S;\ZZ)$ and $\cU^{(k)}$ a universal sheaf over $S \x
  \MM_0^{(k)}$. Let $\mu^{(l_{\z})}(\a)=\mu_-(\a)$ and
  $\mu^{(-1)}(\a)=\mu_+(\a)$. 
  \item Let $z=x^r\a^s \g_1 \cdots \g_a A_1 \cdots A_b$ be any element 
  in $\AA(S)$, where $x\in H_0(S;\ZZ)$ is the generator of the $0$-homology, 
  $\g_i \in H_1(S;\ZZ)$, $\a \in H_2(S;\ZZ)$, $A_i \in H_3(S;\ZZ)$.
  Then we define $\mu^{(k)}(z)$
  as $\mu^{(k)}(x)^r \mu^{(k)}(\a)^s \mu^{(k)}(\g_1) \cdots \mu^{(k)}(\g_a)
  \mu^{(k)}(A_1)\cdots \mu^{(k)}(A_b)$.
\end{itemize}
\end{notation}

Although $\cU^{(k)}$ might not exist, there is always a
well-defined element $p_1({\frak g}_{\cU^{(k)}})$. 
As in~\cite{flips}, we are using the natural complex orientations 
for the moduli spaces. These differ from the natural ones used
in the definition of the Donaldson invariants by a factor
$\f_S(w)=(-1)^{K \cdot w+w^2 \over 2}$. The analogues
of lemma 5.2 and lemma 5.3 of~\cite{flips} are

\begin{lem}
\label{lem:ezlk}
  Let $\g \in H_1(S;\ZZ)$, $\a \in H_2(S;\ZZ)$, $A \in H_3(S;\ZZ)$.
  Put $a=(\z\cdot \a)/2$. Then
  $$ 
    \left\{ \begin{array}{l}
  p_k^*\mu^{(k)}(\a)|_{D_k}=(p_k|_{D_k})^* \left( [\cZ_{\lz-k}]
  \backslash \a + [\cZ_k]
  \backslash \a -a \l_k - c_1(\cF)^2 \backslash \a \right) \sk\\
   p_k^*\mu^{(k)}(\g)|_{D_k}=(p_k|_{D_k})^* \left( [\cZ_{\lz-k}]
  \backslash \g + [\cZ_k]
  \backslash \g -\l_k (c_1(\cF) \backslash \g) \right) \sk\\
   p_k^*\mu^{(k)}(A)|_{D_k}=(p_k|_{D_k})^* \left( [\cZ_{\lz-k}]
  \backslash A + [\cZ_k]
  \backslash A - (\z c_1(\cF)) \backslash A\right) \sk\\
   p_k^*\mu^{(k)}(x)|_{D_k}=(p_k|_{D_k})^* \left( [\cZ_{\lz-k}]
  \backslash x + [\cZ_k]
  \backslash x -{1 \over 4} \l_k^2 \right)
    \end{array} \right. 
  $$ 
\end{lem}

\begin{lem}
\label{lem:pkmu}
  Let $\g \in H_1(S;\ZZ)$, $\a \in H_2(S;\ZZ)$, $A \in H_3(S;\ZZ)$.
  Put $a=(\z\cdot \a)/2$. Then
  $$ 
    \left\{ \begin{array}{l}
   q_{k-1}^*\mu^{(k-1)}(\a) = p_k^*\mu^{(k)}(\a) -aD_k \sk\\
   q_{k-1}^*\mu^{(k-1)}(\g) = p_k^*\mu^{(k)}(\g) - (c_1(\cF)
        \backslash \g) D_k \sk\\
   q_{k-1}^*\mu^{(k-1)}(A) = p_k^*\mu^{(k)}(A) \sk\\
   q_{k-1}^*\mu^{(k-1)}(x) = p_k^*\mu^{(k)}(x) -{1 \over 4}
       (D_k^2 +2\l_k D_k)
    \end{array} \right. 
  $$ 
\end{lem}

We immediately see that it is important to understand the
cohomology classes $e_{\a}=c_1(\cF)^2 \backslash
\a$, $e_{\g}=c_1(\cF) \backslash \g$, and $e_{S}=c_1(\cF)^4 \backslash
[S]$. We write $c_1(\cF)= c_1(F) + \sum \b_i \otimes \b_i^{\#}$, the
K\"unneth decomposition of $c_1(\cF) \in H^2(S\x J)$, where
$\{\b_i\}$ is a basis for $H^1(S)$ and $\{\b_i^{\#}\}$ is the 
dual basis for $H^1(J) \iso H^1(S)^*$. Now we have more explicit 
expressions
  \begin{equation}
    \left\{ \begin{array}{l}
      e_{\a} = -2\sum\limits_{i<j} <\b_i \wedge \b_j, \a>  \b_i^{\#}\wedge
      \b_j^{\#} \in H^2(J) \\
      e_{\g} = \sum <\g, \b_i>  \b_i^{\#} \in H^1(J) \sk\\
      e_{\z A}= (\z c_1(\cF)) \backslash A = \sum <\PD [A] 
       \wedge \b_i, \z>  \b_i^{\#} \in H^1(J) \sk\\
      e_S= \sum\limits_{i,j,k,l}<\b_i \wedge \b_j \wedge \b_k \wedge \b_l, [S]>
      \b_i^{\#}\wedge \b_j^{\#}\wedge \b_k^{\#}\wedge\b_l^{\#} \in H^4(J)   
     \end{array} \right. 
  \label{eqn:the-e's}\end{equation} 
\sk

\begin{thm}
\label{thm:wall-formula}
  Let $\z$ define a wall of type $(w,p_1)$ and
  $d=-p_1-3(1-q)$. Suppose $l_{\z}+ h(\z)+q >0$. For
  $\a \in H_2(S;\ZZ)$, put $a=(\z\cdot \a)/2$. Let
  $z=x^r\a^s \g_1 \cdots \g_a A_1 \cdots A_b \in \AA(S)$ be of degree $2d$.
  Then $\wc(\a)$ is $\f_S(w)$ times
$$
  \sum_{0 \leq k \leq \lz}
  ([\cZ_{\lz-k}]\backslash x+[\cZ_k] \backslash x - {1 \over 4} X^2)^r 
  ([\cZ_{\lz-k}]\backslash \a+[\cZ_k] \backslash \a -e_{\a} +a X)^s
  ([\cZ_{\lz-k}]\backslash \g_1+
$$
$$ +[\cZ_k] \backslash \g_1 +e_{\g_1}X) \cdots 
  ([\cZ_{\lz-k}]\backslash A_b+[\cZ_k] \backslash A_b -e_{\z A_b})  
$$
\sk where $X^N=(-1)^{N-N_{-\z}} \, s_{N-1-N_{\z}-N_{-\z}}
(\cE_{\z}^{\lz-k,k}\oplus (\cE_{-\z}^{k,\lz-k})^{\vee})$,
$s_i(\cdot)$ standing for the Segre class.
\end{thm}

\begin{pf}
  By lemma~\ref{lem:pkmu}, $\mu^{(k-1)}(z)$ is equal to (we omit the pull-backs)
  $$  (\mu^{(k)}(x) -{1 \over 4}
       (D_k^2 +2\l_k D_k))^r (\mu^{(k)}(\a) -aD_k)^s 
    (\mu^{(k)}(\g_1) - e_{\g_1} D_k) \cdots \mu^{(k)}(A_b)
  $$
  which is $\mu^{(k)}(z)$ plus things containing at least one $D_k$.
  So $\mu^{(k-1)}(z) = \mu^{(k)}(z) + D_k \cdot s$, where $s$ is formally
  (recall $\xi_k=-D_k|_{D_k}$)
  $$ {1 \over -\xi_k} \Big( (\mu^{(k)}(x)|_{D_k} -{1 \over 4}
       (\xi_k^2 -2\l_k \xi_k))^r (\mu^{(k)}(\a)|_{D_k} +a\xi_k)^s 
    (\mu^{(k)}(\g_1)|_{D_k} + e_{\g_1} \xi_k) \cdots
  $$
  $$   
  \cdots (\mu^{(k)}(\g_a)|_{D_k} + e_{\g_1} \xi_k) \mu^{(k)}(A_1)|_{D_k} \cdots
    \mu^{(k)}(A_b)|_{D_k} \Big)_0
  $$ 
  where the subindex $0$ means ``forgetting anything not containing at 
  least one $\xi_k$''. So $s$ is (we drop the subindices)
$$
  -{1 \over \xi} \Big( ([\cZ]\backslash x+[\cZ] \backslash x - 
  {1 \over 4}(\xi-\l)^2)^r ([\cZ]\backslash \a+[\cZ] \backslash \a 
  -e_{\a} +a(\xi-\l))^s ([\cZ]\backslash \g_1 +
$$
$$
  + [\cZ] \backslash \g_1
  +e_{\g_1}(\xi-\l)) \cdots ([\cZ]\backslash A_b+[\cZ] \backslash A_b 
  -e_{\z A_b}) \Big)_0
$$
  We need the easy formula (which can be proved by induction)
$$
   {1 \over \xi} \left( (\xi-\l)^N \right)_0={ (\xi-\l)^N
  -(-\l)^N  \over \xi}= \sum_{i=0}^{N-1} (-\l)^i (\xi-\l)^{N-i-1}
$$
  As $\xi-\l$ is the tautological bundle corresponding to
  $\cE_{-\z}^{k,\lz-k}$ (see items 5 to 7 in 
  notation~\ref{not:wall}), we have 
$$    
    \left\{ \begin{array}{l} \l^u=s_{u-N_{\z}}(\cE_{\z}^{\lz-k,k}) \cdot
   \l^{N_{\z}} + O(\l^{N_{\z}-1}) \\
   (\xi-\l)^u=s_{u-N_{-\z}}(\cE_{-\z}^{k,\lz-k}) \cdot
   (\xi-\l)^{N_{-\z}} + O((\xi-\l)^{N_{-\z}-1}) 
    \end{array} \right. 
$$
Evaluating (and doing the sum from $k=0$ to
$k=\lz$) we get the statement of the theorem where
\begin{eqnarray*}
  X^N &=&-\sum (-1)^i \, s_{i-N_{\z}}(\cE_{\z}^{\lz-k,k}) \cdot
  s_{N-i-1-N_{-\z}}(\cE_{-\z}^{k,\lz-k}) =\\
   &=& \sum (-1)^{N-N_{-\z}} \, s_{i-N_{\z}}(\cE_{\z}^{\lz-k,k}) \cdot
  s_{N-i-1-N_{-\z}}((\cE_{-\z}^{k,\lz-k})^{\vee})= \\
  &=& (-1)^{N-N_{-\z}}\, 
  s_{N-1-N_{\z}-N_{-\z}}(\cE_{\z}^{\lz-k,k}\oplus
  (\cE_{-\z}^{k,\lz-k})^{\vee}).
\end{eqnarray*}
\end{pf}

An immediate corollary which generalises~\cite[theorem 5.4]{flips} is

\begin{cor}
\label{cor:wall-formula}
  Let $\z$ define a wall of type $(w,p_1)$ and
  $d=-p_1-3(1-q)$. Suppose $l_{\z}+ h(\z)+q >0$. For
  $\a \in H_2(S;\ZZ)$, put $a=(\z\cdot \a)/2$. Then 
  $\mu_+(\a^d)-\mu_-(\a^d)$ is equal to
  $$
  \sum (-1)^{h(\z)+\lz+j} \, {d! \over j! b! (d-j-b)!} a^{d-j-b}
  ([\cZ_{\lz-k}]\backslash \a+[\cZ_k]\backslash \a)^j \cdot
  e_{\a}^b \cdot
  s_{2\lz-j +q-b}(\cE_{\z}^{\lz-k,k}\oplus
  (\cE_{-\z}^{k,\lz-k})^{\vee}),
  $$
  where the sum runs through $0 \leq j \leq 2\lz$, $0 \leq b
  \leq q$, $0 \leq k \leq \lz$. As $\mu_+(\a^d)-\mu_-(\a^d)$ 
  is computed using the complex orientation, we have that
  $\wc(\a^d)=\f_S(w)(\mu_+(\a^d)-\mu_-(\a^d))$.
\end{cor}

\begin{rem}
\label{rem:ko}
 In Kotschick notation~\cite{Kotschick1}, $\e(\z,w)=(-1)^{({\z-w \over 2})^2}$. So
 $\f_S(w) (-1)^{h(\z)}= (-1)^{d+q}\e(\z,w)$.
\end{rem}

\begin{thm}
\label{thm:lz=0}
  Let $\z$ define a wall of type $(w,p_1)$ and
  $d=-p_1-3(1-q)$. Suppose 
  $\lz+h(\z)+q=0$ i.e. $l_{\z}=0$ and $h(\z)+q=0$. For
  $\a \in H_2(S;\ZZ)$, put $a=(\z\cdot \a)/2$. Let
  $z=x^r\a^s \g_1 \cdots \g_a A_1 \cdots A_b \in \AA(S)$ be of degree $2d$.
  Then $\wc(\a)$ is $\f_S(w)$ times
$$
  (- {1 \over 4} X^2)^r (-e_{\a} +a X)^s (e_{\g_1}X) \cdots (-e_{\z A_b})  
$$
where $X^N=s_{N-N_{-\z}}(\cE_{-\z}^{0,0}) = (-1)^{N-N_{-\z}}\,
s_{N-1-N_{\z}-N_{-\z}}
(\cE_{\z}^{0,0}\oplus (\cE_{-\z}^{0,0})^{\vee})$.
\end{thm}

\begin{pf}
  Now $\MM_+$ is
  $\MM_-$ with an additional connected component $D=E^{0,0}_{-\z}$ which
  is a $\PP^{d-q}$-bundle over $J$, since $E^{0,0}_{\z} = \emptyset$.
  The universal bundle over $E^{0,0}_{-\z}$ is given by an extension
  $$
    \exseq{\pi^*\cO_{S \x J}(\pi_1^* L-\cF) \otimes p^* \l}
    {\cU}{\pi^*\cO_{S \x J}(\cF)},
  $$
  where $\pi: S \x E^{0,0}_{-\z} \ar S \x J$ and $p:S \x E^{0,0}_{-\z}
  \ar E^{0,0}_{-\z}$ are projections and $\l$ is the tautological line
  bundle. From this 
  $$ 
    \left\{ \begin{array}{l}
  \mu (\a)|_D= a\l - e_{\a}  \sk\\
  \mu (\g)|_D= \l e_{\g} \sk\\
  \mu (A)|_D= -e_{\z A}\sk\\
  \mu (x)|_D= -{1 \over 4} \l^2
    \end{array} \right. 
  $$ 
  with notations as in theorem~\ref{thm:wall-formula}. As in the proof of
  theorem~\ref{thm:wall-formula}, $\l^u=s_{u-N_{-\z}}(\cE_{-\z}^{0,0}) \cdot
   \l^{-N_{\z}} + O(\l^{-N_{\z}-1})$, so the expression of the statement 
  of the theorem follows with 
  $X^N=s_{N-N_{-\z}}(\cE_{-\z}^{0,0})$.
\end{pf}

  The next step is to find more handy expressions for the set of classes given
  by~\eqref{eqn:the-e's}.

\begin{lem}
\label{lem:wedge}
  Let $S$ be a manifold with $b^+=1$. Then there is a (rational)
  cohomology class   $\S \in H^2(S)$ such that the image of $\wedge:
  H^1(S) \otimes H^1(S) \ar H^2(S)$ is $\QQ[\S]$. Also $e_S=0$.
\end{lem}

\begin{pf}
  Let $\b_1,\b_2 ,\b_3 ,\b_4  \in H^1(S)$. If $\b_1 \wedge \b_2 \wedge
  \b_3 \wedge \b_4 \ne 0$ then the image of $\wedge:H^1(S) \otimes
  H^1(S) \ar H^2(S)$ contains the subspace $V$ generated by $\b_i \wedge
  \b_j$, which has dimension $6$, with $b^+=3$ and $b^-=3$. This is absurd,
  so $\b_1 \wedge \b_2 \wedge
  \b_3 \wedge \b_4 = 0$. Then $e_S=0$.

  Now let $\S_1 =\b_1 \wedge \b_2$, $\S_2 =\b_3 \wedge \b_4 \in
  H^2(S)$. Then $\S_1^2=\S_2^2=0$ together with the fact that $b^+=1$
  imply that $\S_1 \cdot \S_2 \neq 0$ unless $\S_1$ and $\S_2$ are
  proportional. Since $\S_1 \cdot \S_2 = 0$ by the above, it must be
  the case that $\S_1$ and $\S_2$ are proportional.
\end{pf}

\begin{rem}
  If $S \to C_g$ is a ruled surface with $q>0$ and fiber class $f$, then $\S=f$.
  Note also that the class $\S$ does not change under blow-ups.
\end{rem}

Now write $\seq{\b}{1}{2q}$ for a basis of $H^1(S)$ and fix a generator
$\S$ of the image of $\wedge: H^1(S) \otimes H^1(S) \ar H^2(S)$. 
Let $\seq{\d}{1}{2q}$ be the dual basis for $H_1(S)$.
Put $\b_i \wedge \b_j=a_{ij}
\S$. The Jacobian of $S$ is $J= H^1(S;\RR)/H^1(S;\ZZ)$, so
naturally $H^1(J) \isom H^1(S)^*$. Let $\cL \ar S
\x J$ be the universal bundle 
parametrising divisors homologically equivalent to zero. Then
$E=c_1(\cL)=\sum \b_i \otimes \b_i^{\#}$, with $\b_i^{\#}$ corresponding
to $\d_i$ under the isomorphism $H^1(J) \iso H_1(S)$. So
  \begin{equation} 
    \left\{ \begin{array}{l}
    e_{\a}=-2 \sum\limits_{i<j} a_{ij} (\S \cdot \a) \b_i^{\#}\wedge
    \b_j^{\#}=- 2 ( \S \cdot \a) \o \\
    e_{\d_i} =\b_i^{\#} \sk\\
    e_{\z \b_i} = \sum (\S\cdot \z) a_{ij} \b_j^{\#} = (\S\cdot\z)
    i_{\b_i} \o
    \end{array} \right. 
  \label{eqn:the-e's2}\end{equation}
where we write $\o = \sum\limits_{i<j} a_{ij}(\b_i^{\#}\wedge
 \b_j^{\#}) \in H^2(J)$, which is an element 
independent of the chosen basis. We also have implicitly assumed 
$H_3(S) \iso H^1(S)$ through Poincar\'e duality, in the third line. 
We define
$$
  F: \AA(S) \ar \L^*H_1(S) \otimes \L^*H_3(S) \ar \QQ
$$
given by projection followed by the map sending
$\g_1 \wedge \cdots \wedge \g_a \otimes A_1 \wedge \cdots
  \wedge A_b$ to zero when $a+b$ is odd and to 
$$
 \int_J (\g_1 \wedge \cdots \wedge \g_a 
\wedge i_{A_1} \o \wedge \cdots
\wedge i_{A_b} \o \wedge \o^{q-(a+b)/2})
$$
when $a+b$ is even. We note that we always can find a basis $\seq{\b}{1}{2q}$ with 
$$
\o=a_1 \b_1^{\#}\wedge \b_2^{\#} + a_2 \b_3^{\#}\wedge \b_4^{\#} + \cdots
a_r \b_{2r-1}^{\#}\wedge \b_{2r}^{\#},
$$
where $a_i \neq 0$ are integers and $r \leq q$.
So if $\o$ is degenerate, $F(1) =\int_J \o^q=0$. In general, for a basis element
$z=x^r\a^s\d_{i_1}\cdots \d_{i_a}\b_{j_1} \cdots \b_{j_b}$, $F(z)$ is 
zero unless $z$ contains $\d_{2r+1}\cdots \d_{2q}$, and for every pair
$(2i-1,2i)$, $1 \leq i \leq r$, either
$\d_{2i-1}\d_{2i}$, $\b_{2i-1}\b_{2i}$, 
$\d_{2i-1}\b_{2i-1}$, $\d_{2i}\b_{2i}$ or nothing.
In any case, for subsequent use, we set
$$
  \vol = { 1\over q!} \int_J \o^q.
$$
The number $\vol$ depends on the choice of $\S$, as when $\S$ is changed to $r\S$,
$\vol$ is changed to $r^{-q}\vol$. The final expressions we get for the 
wall-crossing terms are (as expected) independent of this choice.
Also we are going to need the following

\begin{prop}
\label{prop:segre}
  For any sheaf $\cF$ on any complex variety, 
  the Segre classes of $\cF$ are given by $s_t(\cF)=c_t(\cF)^{-1}$. For
  the relationship between the Chern classes of $\cF$ and its Chern
  character, write $a_i$ for $i!$ times 
  the $i$-th component of $\ch \cF$. Then
  $$
    c_n(\cF)={1 \over n!}\left|{\begin{array}{ccccc} a_1 & n-1 & 0 &
    \cdots & 0 \\ a_2 & a_1 & n-2 & \cdots & 0 \\ \vdots & & \ddots & &
    \vdots \\ \vdots & & & \ddots & 1 \\ a_n & a_{n-1} &a_{n-2} & \cdots &
    a_1 \end{array}}\right|
  $$
  and
  $$
    s_n(\cF)={1 \over n!}\left|{\begin{array}{ccccc} -a_1 & -(n-1) & 0 &
    \cdots & 0 \\ a_2 & -a_1 & -(n-2) & \cdots & 0 \\ \vdots & & \ddots & &
    \vdots \\ \vdots & & & \ddots & 1 \\ (-1)^na_n & (-1)^{n-1}a_{n-1}
     &(-1)^{n-2}a_{n-2} & \cdots &
    -a_1 \end{array}}\right|
  $$
\end{prop}

\section{The case $l_{\z}=0$}
\label{sec:lz=0}

In this section we are going to compute $\wc$ in the case $\lz=0$, i.e. 
when $\z^2 =p_1$. We have the following theorem which 
extends~\cite[theorems 6.1 and 6.2]{flips}~\cite{Kotschick1}.

\begin{thm}
\label{thm:mainl0}
  Let $\z$ be a wall with $\lz=0$. Then $\wc(x^r\a^{d-2r})$ is equal to 
  $$
  \e(\z,w) \sum_{0 \leq b \leq q} (-1)^{r+d} 2^{3q-b-d}{q!  \over
       (q-b)!} {d-2r \choose b}
      ( \z \cdot \a)^{d-2r-b}  ( \S \cdot \a)^b ( \S \cdot \z)^{q-b} \vol,
  $$
  where terms with negative exponent are meant to be zero.
\end{thm}

\begin{pf}
For simplicity of notation, let us do the case $r=0$ (the other case is
very similar).
Recall that $F$ is a divisor such that $2F-L$ is homologically equivalent to $\z$ 
and $J=\Jac^F(S)$ is the Jacobian parametrising divisors homologically equivalent 
to $F$. Then $\cF \subset S \x J$ denotes the universal divisor.
Now $\cE_{\z}=\cE_{\z}^{0,0}=R^1 \pi_*(\cO_{S\x J}(2\cF - \pi_1^*L))$
(with $\pi: S \x J \ar J$ the projection) is a vector bundle over $J$.
We note that $H^0(\cO_S (2 F - L)) =0$ and $H^0(\cO_S (L-2 F)
\otimes K) =0$, as $\z$ is a good wall, so $R^0 \pi_*$ and $R^2 \pi_*$ vanish. Then
$$
  \ch \cE_{\z} = -\ch \pi_! (\cO_{S\x J}(2\cF - \pi^*L)) = 
  -\pi_*(\ch \cO_{S\x J}(2\cF - \pi^*L) \cdot \text{Todd }T_S) =
$$
\begin{equation}
  = - ({\z^2 \over 2} -{\z\cdot K \over 2} +1 -q) + 
 e_{K-2\z} -{2 \over 3} e_S = \rk(\cE_{\z}) + e_{K-2\z},
\label{eqn:rar}\end{equation}
since $e_S=0$ (lemma~\ref{lem:wedge}). 
A fortiori $\ch \cE^{\vee}_{-\z} =- ({\z^2 \over 2}
+{\z\cdot K \over 2} +1 -q) - e_{K+2\z}$ and
$$
 \ch( \cE_{\z} \oplus \cE_{-\z}^{\vee}) = (-\z^2 +2q-2)- 4e_{\z}.
$$
From proposition~\ref{prop:segre}, 
$s_i(\cE_{\z} \oplus \cE_{-\z}^{\vee})= {4^i \over i!}
e_{\z}^i$.
This together with theorem~\ref{thm:lz=0} or theorem~\ref{thm:wall-formula} 
(depending on whether $h(\z)+q$ is zero or not) gives
\begin{eqnarray*}
\wc(\a^d) &=& \f_S(w)\sum_{0 \leq b \leq q} (-1)^{h(\z)}
         {d \choose b} a^{d-b} e_{\a}^b \cdot
          s_{q-b}(\cE_{\z}\oplus \cE_{-\z}^{\vee})= \\
   &=& \f_S(w)\sum_{0 \leq b \leq q} (-1)^{h(\z)} {d \choose b} a^{d-b}
       e_{\a}^b \cdot {4^{q-b} \over (q-b)!} e_{\z}^{q-b}= \\
   &=&  \f_S(w)\sum_{0 \leq b \leq q} (-1)^{h(\z)+q} {2^{3q-b-d} \over
       (q-b)!} {d \choose b}
      ( \z \cdot \a)^{d-b}  ( \S \cdot \a)^b ( \S \cdot \z)^{q-b} \o^q,
\end{eqnarray*}
using~\eqref{eqn:the-e's2}. Now we substitute
$\f_S(w) (-1)^{h(\z)}= (-1)^{d+q}\e(\z,w)$ (remark~\ref{rem:ko}), to 
get the desired result. 
\end{pf}

We can also generalise introducing classes of odd degree. If
$z=x^r\a^s\g_1 \cdots \g_a A_1 \cdots A_b$ with $d=s+2r+{3 \over 2}a 
+{1 \over 2}b$, then theorem~\ref{thm:lz=0} or theorem~\ref{thm:wall-formula} gives
\begin{eqnarray*}
\wc(z) &=&  \f_S(w)(- {1 \over 4} X^2)^r (-e_{\a} +a X)^s 
(e_{\g_1}X) \cdots (-e_{\z A_b}) \\
&=&  \f_S(w) 
  \sum_{j} (- {1 \over 4})^r (-1)^b 
  {s \choose j} a^{s-j} 2^j ( \S \cdot \a)^j
 (\S\cdot\z)^b  \o^j \, \g_1 \cdots \g_a \,
 i_{A_1}\o \cdots i_{A_b} \o \cdot \\ & &
 \cdot (-1)^{2r+a+s-j-N_{-\z}} \, s_{2r+a+s-j-1-N_{\z}-N_{-\z}} 
 (\cE_{\z} \oplus \cE_{-\z}^{\vee}).  
\end{eqnarray*}
Now $2r+a+s-1-N_{\z}-N_{-\z}=q-(a+b)/2$, so
\begin{eqnarray*}
\wc(z) &=&  \e(\z,w) \sum_{j} (-1)^{r+d+b} \, 2^{3q-d-b-j} \, 
  {s \choose j} {F(z) \over (q-(a+b)/2-j)!} \cdot \\
   & & \cdot(\z\cdot\a)^{s-j} ( \S \cdot \a)^j (\S\cdot\z)^{q+(b-a)/2-j}\vol.
\end{eqnarray*}

\section{The case $l_{\z}=1$}
\label{sec:lz=1}

Now we want to compute $\wc$ in the case $\lz=1$, i.e. 
when $\z^2 =p_1 + 4$. In this case, $J \x H_{\lz-k}\x  H_k \iso J \x S$,
both for $k=0$ and $k=1$. 
The universal divisor $\cZ_1 \subset S \x H_1 =S \x S$ is
the diagonal $\D$.
Let again $\cL \to S \x J$ be the universal bundle parametrising divisors 
homologically equivalent to zero, so $\cF= \pi_1^* F + \cL$.
With this understood, we have the following 
easy extension of~\cite[lemma 5.11]{flips}

\begin{lem}
\label{lem:5.11}
  Let $\Hom=\Hom(I_{\cZ_k},I_{\cZ_{\lz-k}})$ and
  $\Ext^1=\Ext^1(I_{\cZ_k},I_{\cZ_{\lz-k}})$, $\pi_1$, $p$ and $\pi_2$ be the 
  projections from $S \x (J \x H_{\lz-k}\x H_k)$ to $S$, $S \x J$ and $J \x 
  H_{\lz-k}\x  H_k$, respectively. Let $E=c_1(\cL)$. Then we have
  the following exact sequences
  $$ 
    \exseq{R^1\pi_{2*}(p^*(\z+2E)\otimes
   \Hom)}{\cE_{\z}^{\lz-k,k}}{\pi_{2*}(p^*(\z+2E)\otimes \Ext^1)}
  $$
  $$ 
    \exseq{\pi_{2*}(p^*(\z+2E)\otimes
   \cO_{\cZ_{\lz-k}})}{R^1\pi_{2*}(p^*(\z+2E)\otimes
   \Hom)}{R^1\pi_{2*}(p^*(\z+2E))}
  $$
where the last sheaf is $M_{\z}=R^1 \pi_{2*}(\cO_{S\x J}(2\cF - \pi_1^*L))$,
which is a line bundle over
$J$ with $ch M_{\z} = rk M_{\z} + e_{K-2\z}$ (computed in equation~\eqref{eqn:rar}).
\end{lem}

We apply this lemma to our case $\lz=1$. Then
for $k=0$, $\Hom=\Hom(\cO,I_{\D})=I_{\D}$, $\Ext^1=\Ext^1 (\cO,I_{\D})=0$, 
and for $k=1$, 
$\Hom=\Hom(I_{\D},\cO)=\cO_{S\x S}$, $\Ext^1=\Ext^1(I_{\D},\cO)=\cO_{\D}(\D)$.
Using lemma~\ref{lem:5.11} and the fact 
$\pi_{2*}(\cO_{\D}(\D))=\cO_S(-K)$, we get
$$\ch \cE^{1,0}_{\z}=\ch M_{\z} + \ch \z \, \ch 2E$$
$$\ch \cE^{0,1}_{\z}=\ch M_{\z} + \ch (\z-K)\,  \ch 2E$$
We recall from notation~\ref{not:wall},
\begin{eqnarray*}
 \cE_{\z}^{\lz-k,k}&=& \cE{\text{xt}}_{\pi_2}^1( \cO_{S \x (J\x H_1 \x H_2)} (\pi_1^* L -\cF) \otimes I_{\cZ_2},
\cO_{S \x (J\x H_1 \x H_2)} (\cF) \otimes I_{\cZ_1}) \\
 &=&  \cE{\text{xt}}_{\pi_2}^1( I_{\cZ_2},
\cO_{S \x (J\x H_1 \x H_2)} (\z + 2E) \otimes I_{\cZ_1}), \\
\cE_{-\z}^{k,\lz-k} &=& \cE{\text{xt}}_{\pi_2}^1(I_{\cZ_1},
  \cO_{S \x (J\x H_1 \x H_2)} (-\z - 2E) \otimes I_{\cZ_2}).
\end{eqnarray*}
Then we have $\ch \cE^{1,0}_{-\z}=\ch M_{-\z} + \ch (-\z) \, \ch (-2E)$
and $\ch \cE^{0,1}_{-\z}=\ch M_{-\z} + \ch (-\z-K)\,  \ch (-2E)$. So
$$
\ch (\cE^{1,0}_{\z} \oplus (\cE^{0,1}_{-\z})^{\vee}) =
(-\z^2+2q-2) -4e_{\z} +2\ch \z \ch 2E + {K^2\over 2} +K\z +K(1+2E+2E^2)
$$
$$
  \ch (\cE^{0,1}_{\z} \oplus (\cE^{1,0}_{-\z})^{\vee}) =
  (-\z^2+2q-2) -4e_{\z} +2\ch \z \ch 2E + {K^2\over 2} -K\z - K(1+2E+2E^2)
$$
\begin{equation}
\label{eqn:what}
\end{equation}

We shall compute $s_i=s_i (\cE^{1,0}_{\z} \oplus (\cE^{0,1}_{-\z})^{\vee}) +
s_i (\cE^{0,1}_{\z} \oplus (\cE^{1,0}_{-\z})^{\vee})$, as a class on
$J \x S$.  
From proposition~\ref{prop:segre}, $s_i$ is an
polynomial expression on  $a_i^{(k)}=i! \, \ch_i (\cE^{1-k,k}_{\z} \oplus (
\cE^{k,1-k}_{-\z})^{\vee})$, $k=0,1$. Furthermore, $s_i$ is invariant
under $K \mapsto -K$, and hence an even 
function of $K\z$, $K$, $KE$ and $KE^2$.  
Now the only non-zero even combinations of $K\z$, $K$, $KE$ and
$KE^2$ are $1$ and $K \cdot K$. The first consequence is that we
can ignore $K\z$, $KE$ and
$KE^2$ in $a_i^{(k)}$ for the purposes of computing $s_i$.
So we can suppose
$$
   \left\{ \begin{array}{l}
    a_1^{(k)} = -4e_{\z}+2\z +4E+ (-1)^k K \sk \\
    a_2= 2\z^2 +8E^2 +K^2 +8E\z \sk \\ 
    a_3= 24 E^2 \z =24 e_{\z}[S]
    \end{array} \right. 
$$ 
where $a_i=a_i^{(0)}=a_i^{(1)}$  for $i \geq 2$, and $a_i=0$ for $i \geq 4$
(here we have used that $E^3=0$ and $E^4=0$ as a consequence
of lemma~\ref{lem:wedge}).
Put $a_1 = -4e_{\z}+2\z+4E$ and define
  $$
    I_n = \left|{\begin{array}{ccccc} -a_1 & -(n-1) &
    \cdots & 0 \\ a_2 & -a_1 & \cdots & 0 \\ -a_3 & a_2 & \cdots & 0 \\ 
    0 & -a_3 & \cdots & 0 \\ \vdots & & \ddots &
    \vdots \\ 0 & 0&  \cdots & -a_1 \end{array}}\right|
  $$
and $I_n^{(k)}$ defined similarly with $a_1^{(k)}$ in the place of $a_1$.
Then by proposition~\ref{prop:segre}, $n !\,s_n= I_n^{(0)}+I_n^{(1)}$.
Easily we have $n! \, s_n=2I_n +2 {n \choose 2}K^2 (4e_{\z})^{n-2}$.
Now we can look for an inductive formula for $I_n$. For $n \geq 2$,
\begin{eqnarray*}
   I_n &=& -a_1 I_{n-1} + (n-1) (2\z^2 +K^2 +8E\z) (4e_{\z}-4E)^{n-2}
    +(n-1)8 E^2 I_{n-2} \\ & & -6(n-1)(n-2)[S](4e_{\z})^{n-2} = \\
    &=& -a_1 I_{n-1} + (n-1) (2\z^2 +K^2 +8E\z) (4e_{\z})^{n-2}
     -(n-1)8E\z (n-2)(4E)(4e_{\z})^{n-3} \\ 
     & & +(n-1)8 E^2 I_{n-2} -6(n-1)(n-2)[S](4e_{\z})^{n-2} = \\
    &=& -a_1 I_{n-1} + (n-1) (2\z^2 +K^2 +8E\z) (4e_{\z})^{n-2}
    -8(n-1)(n-2)[S](4e_{\z})^{n-2} \\
     & & +(n-1)8 E^2 \big( (4e_{\z})^{n-2} -(n-2)(4e_{\z})^{n-3}2\z \big) 
    -6(n-1)(n-2)[S](4e_{\z})^{n-2} = \\
    & =& -a_1 I_{n-1} + (n-1) (4e_{\z})^{n-2}
   (2\z^2 +K^2 +8E\z +8E^2-18(n-2)[S]).
\end{eqnarray*}
In the first equality we have used that $I_n P =(4e_{\z}-4E)^n P$, for any 
$P \in H^i(S) \otimes H^j(J)$ with $i=3,4$. In the third equality we use
that $I_n E^2 =(4e_{\z}-2\z)^n E^2$. With this inductive formula for $I_n$,
we get, for $n \geq 2$, 
$$
   I_n = (4e_{\z}-2\z-4E)^n +\sum_{i=2}^n (4e_{\z}-2\z-4E)^{n-i}
    (i-1)(4e_{\z})^{i-2} 
    \Big( 2\z^2 +K^2 
$$
\begin{equation} +8E\z +8E^2-18(i-2)[S] \Big).
\label{eqn:lamia}
\end{equation}
Now for any $k \geq 0$ (we always 
understand ${k \choose i}=0$ if either $i<0$ or $i>k$),
\begin{eqnarray*}
(4e_{\z}-2\z -4E)^k &=& (4e_{\z})^k +k (-2\z -4E)(4e_{\z})^{k-1} +
{k \choose 2} (-2\z -4E)^2(4e_{\z})^{k-2} \\ & &
+ {k \choose 3}3(-2\z)(-4E)^2 (4e_{\z})^{k-3} = \\
& = & (4e_{\z})^k +k (-2\z -4E)(4e_{\z})^{k-1} \\ & & +
{k \choose 2} (4\z^2 +16E\z +16E^2) (4e_{\z})^{k-2} 
-24  {k \choose 3} [S](4e_{\z})^{k-2}.
\end{eqnarray*}
Substituting this into~\eqref{eqn:lamia}, we have
\begin{eqnarray*} 
   I_n &=& (4e_{\z})^n +n (-2\z -4E)(4e_{\z})^{n-1} +
   {n \choose 2} (4\z^2 +16E\z +16E^2) (4e_{\z})^{n-2} \\ & &
   -24  {n \choose 3} [S](4e_{\z})^{n-2}
   +\sum_{i=2}^n \Big( (4e_{\z})^{n-i}
    (i-1)(4e_{\z})^{i-2} (2\z^2 +K^2 -18(i-2)[S]) \\ & &
   + \big( (4e_{\z})^{n-i} + (n-i)(-4E)(4e_{\z})^{n-i-1} \big)
    (i-1)(4e_{\z})^{i-2} 8E\z \\ & &
   + \big( (4e_{\z})^{n-i} + (n-i)(-2\z)(4e_{\z})^{n-i-1} \big)
    (i-1)(4e_{\z})^{i-2} 8E^2 \Big) =\\
   &=& (4e_{\z})^n -n (2\z +4E)(4e_{\z})^{n-1} + (4e_{\z})^{n-2} \Big[ 
   {n \choose 2} (4\z^2 +16E\z +16E^2) 
   -24  {n \choose 3} [S] \\ & &
   +\sum_{i=2}^n (i-1) \big( 2\z^2 +K^2 -18(i-2)[S] + 8E\z +8E^2 - 
   8(n-i)[S] -4(n-i)[S] \big) \Big].
\end{eqnarray*}
Putting this into the expression for $s_n$, we get
\begin{eqnarray*}
   s_n &=& {2 \over n!} \Big( (4e_{\z})^n -n(2\z+4E)(4e_{\z})^{n-1}
   +(4e_{\z})^{n-2} \big[ {n \choose 2}  (6\z^2 + 24E\z +24 E^2 +K^2)
   \\
   & & - 24 {n \choose 3}[S] +\sum_{i=2}^n (i-1) (36-12n-6i)[S] \big] \Big) 
   +{2 \over n!} {n \choose 2}K^2 (4e_{\z})^{n-2}.
\end{eqnarray*}
The expression in the summatory adds up to $-48{n \choose 3}$, so finally 
$$ 
  s_n =2{(4e_{\z})^n \over n!}  -(4\z+8E) {(4e_{\z})^{n-1} \over (n-1)!} +
  (6\z^2 + 2 K^2 + 24E\z+24E^2 ) {(4e_{\z})^{n-2} \over (n-2)!} 
$$
\begin{equation}
  - 24[S]{(4e_{\z})^{n-2} \over (n-3)!},
\label{eqn:sn}
\end{equation}
for $n \geq 2$ (where the last summand is understood to be zero when $n=2$).
This expression is actually valid for $n \geq 0$ under the proviso that the
terms with negative exponent are zero.

\begin{thm}
\label{thm:mainl1}
  Let $\z$ be a wall with $\lz=1$. 
  Then $\wc(x^r\a^{d-2r})$ is equal to  $\e(\z,w)$ times 
  $$  \sum_{b=0}^q 
     (-1)^{r+d+1} 2^{3q-b-d} \Big[   
      ( \z \cdot \a)^{d-2r-b} \Big( {d-2r \choose b}(6\z^2+2K^2-24q-8r)
       + 8{d-2r \choose b+1}{b+1 \choose 1} \Big) + 
  $$
  $$ +  8 ( \z \cdot \a)^{d-2r-b-2}\a^2 {d-2r \choose b+2}{b+2 \choose 2}
    \Big]
    ( \S \cdot \a)^b ( \S \cdot \z)^{q-b} {q!\over
       (q-b)!} \vol,
  $$
  where terms with negative exponent are meant to be zero.
\end{thm}

\begin{pf}
  By theorem~\ref{thm:wall-formula}, $\wc(x^r\a^{d-2r})= \f_S(w)
  ([S]-{1 \over 4}X^2)^r(\a -e_{\a}+aX)^{d-2r}$ evaluated on $J\x S$,
  where 
\begin{eqnarray*}
   X^N &=&(-1)^{N-N_{-\z}} \, \left( s_{N-1-N_{\z}-N_{-\z}}
 (\cE^{1,0}_{\z} \oplus (\cE^{0,1}_{-\z})^{\vee}) +
  s_{N-1-N_{\z}-N_{-\z}} (\cE^{0,1}_{\z} \oplus (\cE^{1,0}_{-\z})^{\vee})
  \right) = \\ &=& (-1)^{N-N_{-\z}} \,
  s_{N-1-N_{\z}-N_{-\z}}.
\end{eqnarray*}
   Hence
  $$ 
    \wc(x^r\a^{d-2r})=\f_S(w) \sum_b (-1)^{h(\z)+1} 
    {d-2r \choose b}a^{d-2r-b}(-{1 \over 4})^r e_{\a}^{b-2} 
    \cdot
  $$
  $$ \cdot \Big[ -4r[S]e_{\a}^2 \cdot s_{q-b} +{b \choose 2} \a^2\cdot
   s_{q-b+2}+{b \choose 1} \a (-e_{\a})\cdot s_{q-b+2}+e_{\a}^2 
   \cdot s_{q-b+2} \Big].
  $$
  Substituting the values of $s_n$ from~\eqref{eqn:sn} and using 
  remark~\ref{rem:ko}, we get
\begin{eqnarray*}
\wc(x^r\a^{d-2r}) &=&  \e(\z,w)\sum_b 
     (-1)^{r+d+1} 2^{3q-b-d} \, {d-2r \choose b}
      ( \z \cdot \a)^{d-2r-b} \Big[(6\z^2+2K^2 \\ 
      & & -24q-8r)
      ( \S \cdot \a)^b {( \S \cdot \z)^{q-b}\over
       (q-b)!} + 16 (\z \cdot \a) {b\choose 1}
      ( \S \cdot \a)^{b-1} {( \S \cdot \z)^{q-b+1}\over
       (q-b+1)!} + \\ 
       & & + 32 \a^2 {b \choose 2}
      ( \S \cdot \a)^{b-2} {( \S \cdot \z)^{q-b+2}\over
       (q-b+2)!} \Big]\, \int_J \o^q.
\end{eqnarray*}
Reagrouping the terms we get the desired result.
\end{pf}

This result agrees with theorems 6.4 and 6.5 in~\cite{flips} particularising
for $q=0$ and $r=0,1$. We see from theorem~\ref{thm:mainl1} that the
difference terms $\wc$ do not satisfy in general the simple type 
condition~\cite{KM}.

\begin{rem}
  L.\ G\"ottsche and the author have obtained the same formula of
  theorem~\ref{thm:mainl1} in some examples, like $\CP^1 \x C_1$ 
  ($C_1$ being an elliptic curve) 
  using the simple type condition in limiting chambers. These arguments
  will appear elsewhere.
\end{rem}

\section{General case}
\label{sec:5}

We do not want to enter into more detailed computations of the
wall-crossing formulae, but just to remark that the pattern laid
in~\cite{flips} together with theorem~\ref{thm:wall-formula}
can be used here to obtain partial information of $\wc$. For
instance, we write
$$
  S_{j,b}= \sum_k([\cZ_{\lz-k}]\backslash \a+[\cZ_k]
  \backslash \a)^j \cdot e_{\a}^b \cdot
  s_{2\lz-j +q-b}(\cE_{\z}^{\lz-k,k}\oplus
  (\cE_{-\z}^{k,\lz-k})^{\vee}),
$$
so that corollary~\ref{cor:wall-formula} says
$\wc(\a^d)=\e(\z,w) \sum (-1)^{d+q+\lz+j}{d! \over j!b!(d-j-b)!}a^{d-j-b} S_{j,b}$.
Then we can obtain (compare~\cite[proposition 5.12]{flips})
\begin{eqnarray*}
  S_{2\lz,q} &=&  {(2\lz)! \over \lz!} (\a^2)^{\lz} \, e_{\a}^q \\
  S_{2\lz-1,q} &=&  (-4){(2\lz)! \over \lz!} (\a^2)^{\lz-1} \,
     a \, e_{\a}^q \\
  S_{2\lz,q-1} &=&  4 {(2\lz)! \over \lz!} (\a^2)^{\lz} \,
     e_{\a}^{q-1} \, e_{\z}
\end{eqnarray*}

As an easy consequence of this we get (compare~\cite[theorems 5.13 and 5.14]{flips})

\begin{cor}
  Let $\z$ be a wall of type $(w,p_1)$. Let $\a \in H_2(S;\ZZ)$ and $a=(\z\cdot\a)
  /2$. Then $\wc(\a^d)$ is congruent (modulo 
  $a^{d-2\lz-q+2}$) with
  $$
    \e(\z,w) (-1)^{d+\lz} \, 
    2^q \Big[ a^{d-2\lz-q} { d! \over \lz! (d-2\lz-q)!}
    (\a^2)^{\lz} (\S\cdot\a)^q + 
  $$
   $$
     + 4a^{d-2\lz-q+1} { d! \, q \over \lz! (d-2\lz-q+1)!}
    (\a^2)^{\lz} (\S\cdot\a)^{q-1}(\S\cdot\z) \Big] \vol.
  $$
\end{cor}

\begin{cor}
  In the conditions of the previous corollary, suppose furthermore $d-2r \geq 2\lz+q$.
  Then $\wc(x^r\a^{d-2r})$ is congruent (modulo 
  $a^{d-2r-2\lz-q+2}$) with
  $$
    \f(\z,w)(-1)^{d+\lz+r} \, 2^{q-2r} (-{1 \over 4})^r 
    \Big[ a^{d-2r-2\lz-q} { (d-2r)! \over \lz! (d-2r-2\lz-q)!}
    (\a^2)^{\lz} (\S\cdot\a)^q +
  $$
  $$
     + 4 a^{d-2r-2\lz-q +1} { (d-2r)! \, q \over \lz! (d-2r-2\lz-q+1)!}
    (\a^2)^{\lz} (\S\cdot\a)^{q-1}(\S\cdot\z) \Big] \vol.
  $$
\end{cor}

\section{Conjecture}
\label{sec:conj}

It is natural to propose the following

\noindent {\bf Conjecture.}
  Let $X$ be an oriented compact
  four-manifold with $b^+=1$ and $b_1=2q$ even. Let $w \in
  H^2(X;\ZZ)$. Choose $\S \in H^2(X)$ generating the image of $\wedge:
  H^1(X) \otimes H^1(X) \ar H^2(X)$. Define $\o \in H^2(J)$ such that 
  $e_{\a}=-2(\S\cdot \a)\o$ and put $\vol=\int_J \frac{\o^n}{n!}$. If
  $\z$ defines a wall, then the wall-crossing difference term 
  $\d_{X,\z}^{w,d} (x^r\a^{d-2r})$ only depends on $w$, $d$, $r$,
  $b_1=2q$, $b_2$, $\z^2$, $\a^2$, $(\z \cdot \a)$
  and $(\S\cdot \a)^i(\S\cdot \z)^{q-i} \vol$, $0 \leq i
  \leq q$. The coefficients are universal on $X$.

This is quite a strong conjecture and one can obviously write down
weaker versions. It would allow one to carry out similar arguments to
those in~\cite{Gottsche} and therefore to find out the general shape
of the wall-crossing formulae for arbitrary $X$, involving
modular forms. One should be able to determine then 
all wall-crossing formulae from particular cases.
This and applications to computing the invariants of $\CP^1 \x C_g$ 
($C_g$ the genus $g$ Riemann surface)
will be carried out in following joint work with L.\ G\"ottsche.

\section*{Appendix. Algebraic surfaces with $p_g=0$ and $-K$ effective}

From~\cite{BPV}, the algebraic surfaces with $p_g=0$ and $-K$ effective are
$\CP^2$, ruled surfaces and blow-ups of these.
For the case $q=0$, we have thus $\CP^2$, the Hirzebruch surfaces and their blow-ups.
Not all blow-ups have $-K$ effective, but they are always deformation equivalent to 
one with $-K$ effective.
For the case $q>0$, the minimal models are ruled surfaces over a surface $C_g$ of
genus $g \geq 1$. They have $c_1^2=8(1-g)$. Let $S \to C_g$ be a ruled surface. 
It has $b_2=2$ and $b_1=2g$, so $g=q$. Let $f$ be the class of the
fibre and $\s=\s_{-N}$ the class of the section with negative
self-intersection $\s_{-N}^2=-N \leq 0$. 
Then there is a section $\s_N$ homologically equivalent to $\s_{-N}+Nf$ with square $N$.
Write $X=\PP(V^{\vee})$, for $V \to C_g$ a rank two bundle. 
Then $K=\fra f-2\s$, with $\fra=\s^2 +K_{C_g}$ a divisor on $C_g$ 
(see~\cite[section 5.2]{Hartshorne}).
Therefore $-K$ is effective if and only if $-\fra$ is effective.
The section $\s$ corresponds to a sub-line bundle $L \inc V$ with
$\cO_{C_g}(\s^2)=L^{-2}\otimes \det(V)$. Then $-\fra$ is effective
when $L^2 \otimes \det(V)^{-1} \otimes K_{C_g}^{-1}$ has sections.
We can find examples for any $N$ as long as $N \geq 2(g-1)$.
Again, the non-minimal examples are blow-ups of these, and can be found to have $-K$
effective.

For fixed $q=g >0$, there are only two deformation classes of minimal ruled surfaces,
corresponding to two diffeomorphism types,
the two different $\SS^2$-bundles over $C_g$, one with even
$w_2$, the other with odd $w_2$. 
\begin{itemize}
 \item {\bf $N$ even:} $S$ is diffeomorphic to $S_0=\CP^1 \x C_g$ (and the canonical
 classes correspond). Let $C$ be the 
 homology class of\, $\pt \x C_g$ coming from the diffeomorphism. Then $\s$ is 
 homologous to $C -{N \over 2}f$. The ample cone $C_S$ of $S$ is generated by $f$ and 
 $\s_N= C + {N \over 2}f$ (i.e. it is given by $\RR^+ f + \RR^+ \s_N$).
 Note that the bigger $N$, the smaller the ample cone.
 The wall-crossing terms $\wc$ do not depend on the complex structure of $S$,
 so our results for the case $-K$ effective give the wall-crossing terms for $S_0$ for
 any wall inside $C_S$. Letting $N=2(g-1)$, we actually compute $\wc$ for 
 any $\z= a\, \CP^1 -b \,C$ with $a,b >0$, $a> b (g-1)$ 
 (note that all these walls are good).
 \item {\bf $N$ odd:} $S$ is diffeomorphic to the non-trivial $\SS^2$-bundle over $C_g$.
 Arguing as above, we compute the wall-crossing terms $\wc$ for 
 any $\z= a \, \CP^1 -b \, \s_{-(2g-1)}$ with $a,b >0$, $a> b {2g-1 \over 2}$.
\end{itemize}

\end{document}